\documentclass[11pt,reqno]{amsart}

\usepackage{epsfig}
\usepackage{amssymb,amsmath}
\usepackage{graphicx}

\newcommand{\be}{\begin{eqnarray}}
\newcommand{\ee}{\end{eqnarray}}

\newcommand{\clos}{\mbox{\rm clos}}

\newcommand{\es}{\emptyset}
\newcommand{\half}{\frac{1}{2}}

\newcommand{\R}{{\mathbb R}}

\newcommand{\Z}{{\mathbb Z}}
\newcommand{\Nat}{{\mathbb N}}
\def\N{\Nat}

\newcommand{\disc}{{\mathbb D}}

\def\Nk{{\mathcal N}}
\def\Ok{{\mathcal O}}
\def\Rk{{\mathcal R}}
\def\Diag{{\rm Diag}}
\def\Gam{\Gamma}

\newcommand{\bob}{{\bf b}}

\newcommand{\bx}{{\bf x}}

\newcommand{\lam}{\lambda }
\newcommand{\gam}{\gamma}
\newcommand{\Om}{\Omega}
\newcommand{\Omtil}{\widetilde{\Om}}
\newcommand{\om}{\omega}
\newcommand{\Nktil}{\widetilde{Q}}
\newcommand{\ftil}{\widetilde{f}}
\newcommand{\htil}{\widetilde{h}}
\newcommand{\gatil}{\widetilde{\gam}}
\newcommand{\latil}{\widetilde{\lam}}

\newcommand{\eps}{\varepsilon}

\newcommand{\Bk}{{\mathcal B}}

\newcommand{\Mk}{{\mathcal M}}

\newcommand{\Xk}{{\mathcal X}}

\def\wtil{\widetilde}
\def\Nktil{\wtil{\Nk}}

\newtheorem{theorem}{Theorem}[section]
\newtheorem{lemma}[theorem]{Lemma}
\newtheorem{cor}[theorem]{Corollary}
\newtheorem{prop}[theorem]{Proposition}

\theoremstyle{definition}
\newtheorem{defi}[theorem]{Definition}

\theoremstyle{remark}

\numberwithin{equation}{section}

\begin{document}

\title[Connectedness locus for pairs of affine maps]{Connectedness locus for pairs of affine maps  \\ [1.2ex] and
zeros of power series}

\author{Boris Solomyak }
\address{Boris Solomyak, Box 354350, Department of Mathematics,
University of Washington, Seattle, WA, USA}
\email{solomyak@uw.edu}


\thanks{The author was supported in part by NSF grant DMS-0968879}

\begin{abstract}
We study the connectedness locus $\Nk$ for the family of iterated function systems of pairs of affine-linear maps in the plane (the non-self-similar case).
First results on the set $\Nk$ were obtained in joint work with P. Shmerkin \cite{Pablo-Sol}. 
Here we establish rigorous bounds for the set $\Nk$ based on the study of power series of special form. We also derive some bounds for the 
region of ``$*$-transversality'' which have applications to the computation of Hausdorff measure of the self-affine attractor.  
We prove that a large portion of the set $\Nk$
is connected
and locally connected, and 
conjecture that the entire connectedness locus is connected.
We also prove that the set $\Nk$
has many zero angle ``cusp corners,'' at certain points with algebraic coordinates.
\end{abstract}

\date{\today}


\maketitle

\thispagestyle{empty}

\section{Introduction}

Here we set the notation and discuss earlier results on the set $\Nk$. This section has some overlap with 
the introductory part of \cite{Pablo-Sol}.
Let $E=E(T,\bob)$ be the attractor of the IFS $\{T\bx, T\bx + \bob\}$,
{\em i.e.}, the unique nonempty compact set in $\R^d$ satisfying
\be \label{saf1}
E = TE \cup (TE + \bob).
\ee
Observe that
\be \label{serep}
E(T,\bob) = \Bigl\{\sum_{n=0}^\infty a_n T^n \bob:\ a_n \in \{0,1\} \Bigr\}
\ee
since the right-hand side is well-defined (that is, the sums converge because $T$ is a contraction, and the set is compact and non-empty) and satisfies
(\ref{saf1}).

We can assume that all the eigenvalues of $T$
have spectral (geometric) multiplicity one, and $\bob$ is a cyclic vector
for $T$, that is, $H  := Span\{T^k \bob:\ k\ge 0\} = \R^d$.
There is no loss of generality in making this assumption, since otherwise
we can replace $T$ by the restriction of $T$ to $H$
and consider the corresponding IFS on $H$.

It is well-known (see \cite{hata}) that the set
$E=E(T,\bob)$ is connected if and only if $TE \cap (TE + \bob) \ne \es$.
This easily implies the following criterion for connectedness.
Denote
$$
\Bk = \Bigl\{1 + \sum_{n=1}^\infty b_n z^n:\ b_n \in \{-1,0,1\}\Bigr\}.
$$
The symbol $\disc$ stands for the open unit disk.

\begin{prop}[\cite{Pablo-Sol}] \label{prop-connect} Let $T$ be a linear contraction with
(possibly complex) eigenvalues $\lam_j$, for $j=1,\ldots, m$,
having algebraic 
multipicities $k_j \ge 1$, and geometric multiplicities equal
to one. Let $\bob$ be a cyclic vector for $T$. Then $E(T,\bob)$ is
connected if and only if there exists $f \in \Bk$ such that
\be \label{eq-zero}
f(\lam_j) = \ldots = f^{(k_j-1)}(\lam_j) = 0,\ \ \ j=1,\ldots,m.
\ee
In particular, connectedness does not depend on $\bob$.
\end{prop}

\medskip

From now on, we restrict ourselves to the case $d=2$. Applying an
invertible linear transformation as a conjugacy, we can assume without
loss of generality that $T$ is one of the following:
$$
{\rm (i)}\ \ T = \left[\begin{array}{cc} a & b \\ 
-b & a \end{array} \right],
\ \ \ \ \ {\rm (ii)}\ \ T = \left[\begin{array}{cc} \gam & 0 \\ 0 & \lam 
\end{array} \right],
\ \ \ \ \ {\rm (iii)}\ \ T = \left[\begin{array}{cc} \lam & 1 \\ 0 & \lam
\end{array} \right],
$$
where $\lam, \gam, a, b$ are real, 
$|\lam|, |\gam| < 1$, and $a^2 + b^2 < 1$.
Note that $\lam \ne \gamma$ by the assumption that $T$ has a cyclic
vector.
The following corollary is immediate from Proposition~\ref{prop-connect}.

\begin{cor}[\cite{Pablo-Sol}] \label{cor-connect} Let $E(T,\bob)$ be the attractor of the
IFS $\{T\bx, T\bx + \bob\}$ where $T$ is of the form (i), (ii), or (iii),
and let $\bob$ be a cyclic vector for $T$.

{\bf (a)} In the case (i), the self-affine set $E(T,\bob)$ is connected
if and only if there exists $f \in \Bk$ such that $f(a+ib) = 0$.

{\bf (b)} In the case (ii), the self-affine set $E(T,\bob)$ is connected
if and only if there exists $f \in \Bk$ such that $f(\lam) = f(\gam)=0$.

{\bf (c)} In the case (iii), the self-affine set $E(T,\bob)$ is connected
if and only if there exists $f \in \Bk$ such that $f(\lam)=f'(\lam)=0$.
\end{cor}

Each of the cases leads to a set which we call the {\em connectedness
locus} for the corresponding family of self-affine sets.  Let
$$
\Mk:=\{z = a+ib \in \disc:\ \exists\,f\in \Bk,\ f(z)=0\},
$$
$$
\Nk:=
\{(\gam, \lam) \in (-1,1)^2:\ \exists\,f\in \Bk, \ f(\gam) = f(\lam) =0\},
$$
$$
\Ok:=\{\lam \in (-1,1): \ \exists\,f\in \Bk,\ f(\lam) = f'(\lam)=0 \}.
$$

Thus, $\Mk$, $\Nk$, and $\Ok$ are essentially the sets of parameters for which the attractors in cases (i), (ii), (iii) are connected.
The only difference is that we allow $b=0$ in $\Mk$ and $\gam=\lam$ in $\Nk$ to ensure that these sets are relatively closed in the unit disk.


The set $\Mk$ has been extensively studied as the 
{\em Mandelbrot set for the pair of linear maps}, see e.g.\ \cite{BH,bou,bandt,sxu,solasymp,ERS} and references therein.

Note that in case (i) the attractors are self-similar, which simplifies
some of the considerations. 

This paper
is devoted to the study of the set $\Nk$, or rather, $\Nk \cap (0,1)^2$.
(By symmetry, we can assume that $\lam>0$. However, the case of $\gam<0$
does not reduce to the case of $\lam$ and $\gam$ having the same sign and
we leave it for a future study.)
It is easy to see (\cite{Pablo-Sol}) that
$$
\{(\lam,\gam)\in [0.5,1)^2:\ \lam\gam \ge 0.5\} \subset \Nk\cap (0,1)^2 \subset
[0.5,1)^2.
$$
A picture of the set $\Nk$ is shown in Figure 1 (which also appears in \cite{Pablo-Sol}). It is created by a 
program of Christoph Bandt, similar to the one used in \cite{bandt} to
draw the set $\Mk$. The set $\{(\lam,\gam)\in [0.5,1)^2:\ \lam\gam \ge 0.5\}$
is shaded grey. The algorithm rigorously checks that a point is outside $\Nk$ and paints it ``white.'' The points that are not declared to be ``white'' after a
certain number of iterations are declared to be in $\Nk$ and painted ``black.''
Thus the figure should be viewed as an ``outward approximation'' for $\Nk$.
However, this is not completely accurate; for instance, the apparent 
disconnected pieces of $\Nk$ are a computing artifact, as we show below. 
Another remark is that the computation is very time-consuming near
the diagonal, so the picture is not accurate there.



\begin{figure}[htb]
\centering
\includegraphics[width=0.8\textwidth]{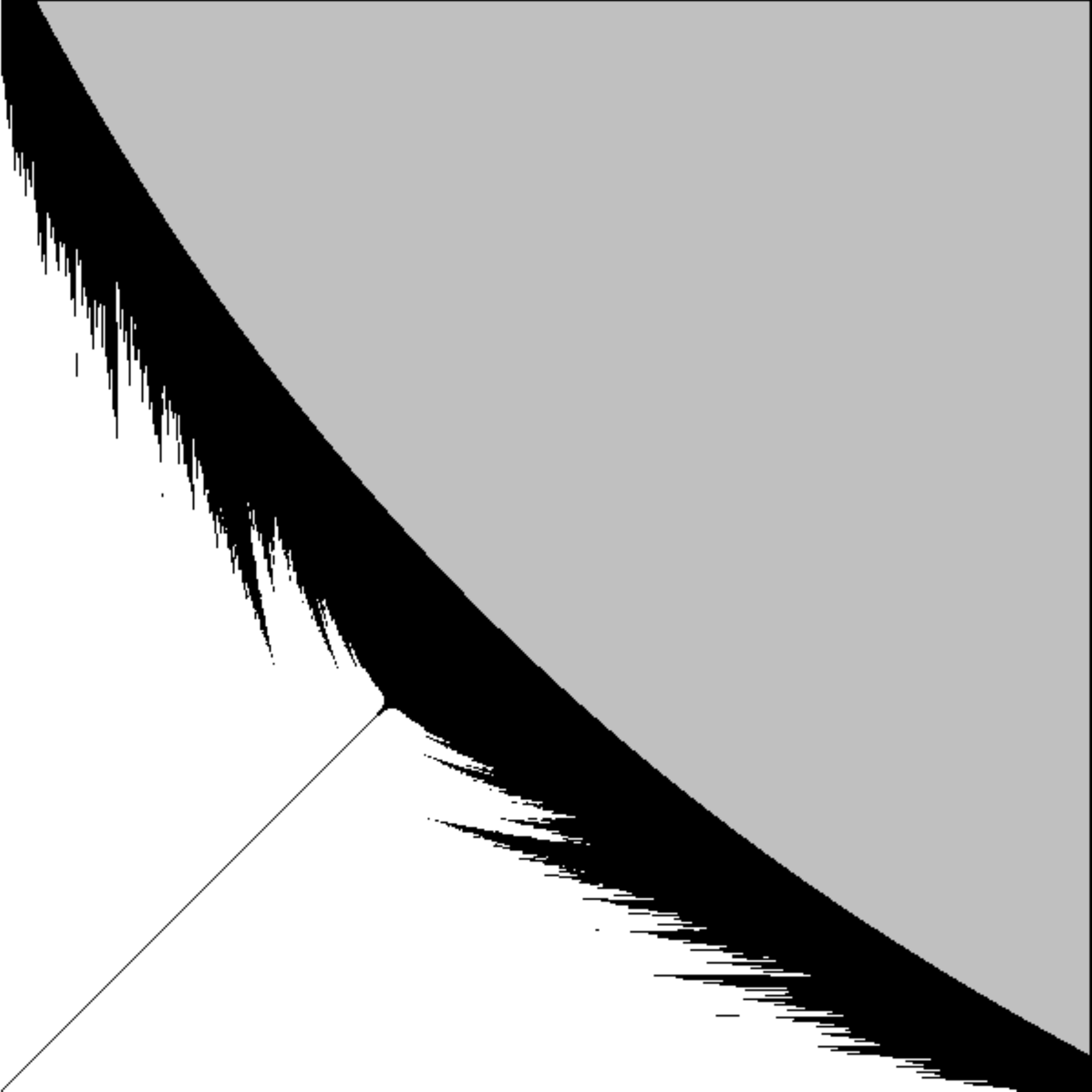}
\caption{Connectedness locus for the family of self-affine sets}
\end{figure}


Next we recall for completeness the results on the set $\Nk$ obtained in \cite{Pablo-Sol}.
Denote $\Diag(F) := \{(\lam,\lam): \,\lam\in F\}$
for $F\subset \R$. We see that the set $\Nk$ has an ``antenna''
$\Gam(\Nk)$, defined in \cite{Pablo-Sol} as the connected component of $\Diag([\half,1))\setminus \clos(\Nk\setminus \Diag(\R))$
containing $(\half,\half)$. In fact, we can consider the set of points on the diagonal which are limit points of $\Nk\setminus\Diag(\R)$. By definition, this
set consists of  those $(\lam,\lam)$ for which there exist $f_n\in \Bk$ with real zeros $\gam_n < \lam_n$ such that $\gam_n\to \lam,\ \lam_n\to \lam$. By compactness of $\Bk$ it follows that there exists $f\in \Bk$ with a double zero at $\lam$, that is, $\lam\in \Ok$. Conversely, a point $(\lam,\lam)$, where $\lam\in (0,2^{-1/2})$
is a double zero of some power series $f\in\Bk$, is in the closure of $\Nk$ if $f$ has infinitely many coefficients not equal to $-1$, since we can
then make an arbitrarily small {\em negative} perturbation of $f$ staying in $\Bk$, which will result in a pair of real zeros close to $\lam$.
(Here we use the fact that $\lam$ is necessarily a local minimum of $f$. It cannot be triple zero, since the smallest triple zero of $f\in \Bk$ is
at least $0.72>2^{-1/2}$ by \cite{BBBP}; see more about this in the next section.)
It follows from \cite{Pablo-Sol} that $\clos(\Nk\setminus \Diag(\R)) \cap \Diag(\R)$ is disconnected, and it is conjectured to have infinitely many connected components.
The ``tip'' of the antenna, that is,  $(\beta,\beta)\in \Gam$ such that $\beta$ is maximal, is found in \cite[Cor.2.10]{Pablo-Sol} with
high accuracy: $\beta =.6684756\pm 10^{-7}$.

Another interesting question concerns the topological structure of $\Nk$. By analogy with the Bandt's conjecture from \cite{bandt},  we expect
that $\Nk\setminus \Diag(\R)$ is contained in the closure of the set of its interior points. It is not obvious even that there exist interior points in the
nontrivial part of $\Nk\setminus \Nk_t$, where $\Nk_t =\{(\gam,\lam) \in (-1,1)^2: \ |\gam||\lam| \ge \half\}$. However, in \cite{Pablo-Sol} it was shown
that a small, but explicitly given, disk around $(2^{-1/2},2^{-1/2})$ is contained in $\Nk$.


\section{Statement of results}

There is another method, which does not involve much computing,
to show that certain regions are disjoint from
$\Nk$. It is based on the idea
that it is much easier to estimate zeros
of power series with ``convex'' restrictions on the coefficients and
uses so-called $(*)$-functions, first introduced in \cite{solerd}. We will
also need their generalizations from \cite{BBBP}.

\begin{defi} \label{mstar}
A power series $h(x) = 1 + \sum_{n=1}^\infty a_n x^n$ is called an
$(m*)$-function if there exist integers $1 \le \ell_1 < \ell_2 < \ldots 
< \ell_m < \infty$ such that $a_{\ell_k}$ are any real numbers for
$k=1,\ldots,m$, and
$$
\begin{array}{ll}
a_n = -1, & 1 \le n \le  \ell_1-1; \\
a_n = (-1)^k, & \ell_{k-1}+1 \le n \le \ell_k-1,\ \ k = 2,\ldots,m; \\
a_n = (-1)^{m+1}, & n \ge \ell_m+1. \end{array}
$$
Moreover, we require that $h$ has exactly $(m+1)$ coefficient sign changes. (It is clear from the assumptions on $a_n$ that the number of
sign changes is at most $(m+1)$, however, it could potentially be less, if for some $j$ we have $\ell_{j+1} = \ell_j+1,\ \ell_{j+2} = \ell_j+2$.)
A $(1*)$-function will be called a $(*)$-function, and a $(2*)$-function
will be called a $(**)$-function.
\end{defi}
Let 
\begin{eqnarray*}
\Nk_+& :=& \Nk\cap \{(\gam,\lam)\in (0,1)^2 :\ \gam < \lam\}\\
          &=& \{(\gam,\lam):\ 0 < \gam < \lam<1 \ \mbox{and there is}\  f \in \Bk,\
f(\gam) = f(\lam) = 0\}.
\end{eqnarray*}
Further, consider
$$
\Bk_{[-1,1]} := 
\Bigl\{1 + \sum_{n=1}^\infty a_n z^n:\ a_n \in [-1,1]\Bigr\}\subset \Bk
$$
and 
$$
\Omega_+ := \{(\gam,\lam):\ 0 < \gam < \lam<1 \ \mbox{and there is}\ f \in \Bk_{[-1,1]},\
f(\gam) = f(\lam) = 0\}.
$$
By definition, $\Nk_+\subset \Omega_+$.
For a power series $f$ with bounded real coefficients, let $\xi_1(f) \le
\xi_2(f) \le \ldots$ denote its positive zeros ordered by magnitude
and counted with multiplicity (for convenience we let $\xi_k(f) = 1$ if
there are fewer than $k$ positive zeros).
In \cite{BBBP} it is proved that for any $k \ge 2$, the smallest
$k$-th order zero $\alpha_k$ of a power series in $\Bk_{[-1,1]}$ is 
algebraic, and the corresponding power series is a $(k*)$-function.
In particular, $\alpha_2 \approx .649138$ is the positive
zero of $2x^5-8x^2+11x-4$.

\thispagestyle{empty}

\begin{prop}\label{prop-om1}
{\bf (a)} The function 
$$
\phi:\ \gam\mapsto \min\{\xi_2(f):\ f \in \Bk_{[-1,1]},\ 
f(\gam) = 0\}
$$
is well-defined on $(.5, \alpha_2)$. It is continuous,
decreasing, and satisfies
$$
\lim_{t\to .5+} \phi(t) = 1,\ \ \ \lim_{t\to \alpha_2-} \phi(t) = \alpha_2.
$$

{\bf (b)} For every $\gam \in (.5,\alpha_2)$ we have
$$
\gam < \lam < \phi(\gam)\ \ \Longrightarrow\ \ (\gam,\lam) \not\in \Om_+.
$$

{\bf (c)} For every $\gam \in (.5,\alpha_2)$ there exists a unique function in $\Bk_{[-1,1]}$ which vanishes at $\gam$ and $\phi(\gam)$.
Moreover, it is a
$(*)$-function
\be \label{eq-star}
h_k^{(a)}(x) = 1 - x - \ldots - x^{k-1} + ax^k + \frac{x^{k+1}}{1-x} \in \Bk_{[-1,1]}
\ee
such  that $h_k^{(a)}(\gam) = h_k^{(a)}(\phi(\gam))=0$. Moreover, $(h_k^{(a)})'(\gam)<0$ and $(h_k^{(a)})'(\phi(\gam))>0$.
\end{prop}

The following table contains some values of the function $\phi$
(rounded-off in such a way that the actual values are slightly
larger):

\medskip

\begin{center}
\begin{tabular}{|c|c|c|c|c|c|c|c|c|c|c|c|c|c|c|}  \hline
$\gam$ & .51 & .52 & .53 & .54 & .55 & .56 & .57 & .58 & .59
\\ \hline
$\phi(\gam)$ & .862 & .831 & .811 & .79 & .77 & .755 & .742 & .728 & .716 
 \\ \hline
\end{tabular}

\smallskip

\begin{tabular}{|c|c|c|c|c|c|c|c|c|c|c|c|c|c|c|}  \hline
$\gam$ & .6 & .61 &
.62 & .63 & .64  \\ \hline
$\phi(\gam)$ & 
.703 & .691 & .68 & .67 & .658  \\ \hline
\end{tabular}

\end{center}


\medskip

Note that the $(*)$-function $h_k^{(a)}$ is in $\Bk$ if and only if $a\in \{-1,0,1\}$. Thus, we get a countable set of points which belong to
$\partial \Om_+\cap \Nk_+$. It turns out that the set $\Nk$ has ``cusp corners'' at these points, as we show in our first main theorem.
This property distinguishes $\Nk$ from  the set $\Mk$, which has spiral points and no corners with
interior angle less than $2\pi/3$ (conjecturally, none at all), see \cite{Sol_locgeom}.

\begin{theorem} \label{th-tips}
Suppose that $\half < \gam_0 < \lam_0 < 1$ and $(\gam_0,\lam_0)$ is such that there is a unique function $h\in \Bk$ which vanishes at
$\gam_0$ and $\lam_0$, and moreover, all the coefficients of $h$ are eventually $+1$ and $h'(\gam_0)< 0,\  h'(\lam_0)> 0$.
Then $z_0 = (\gam_0,\lam_0)$ is a
``tip of a corner'' of the set $\Nk_+$, with zero interior angle.
More precisely, there exist $\delta>0$ and positive constants $C_1$ and $C_2$, 
such that
$$
B_\delta(z_0) \cap \Nk_+ \subset \{(\gam,\lam):\ C_1(\gam_0 - \gam)^\alpha
< \lam-\lam_0 < C_2(\gam_0-\gam)^\alpha\}
$$
where $\alpha = \frac{\log{\lam_0}}{\log{\gam_0}}$. In fact, we can take
$$C_1 = \frac{2|h^{\prime}(\gam_0)|^{\alpha}(1-\gam_0)^{\alpha}}{2^{\alpha}3h^{\prime}(\lam_0)} \ \ \ \mbox{and}\ \ \ 
C_2 = \frac{3^{\alpha}2|h^{\prime}(\gam_0)|^{\alpha}}{2^{\alpha}(1-\lam_0)h^{\prime}(\lam_0)}.
$$
In particular, these conditions are satisfied if $h$
is a $(*)$-function, {\em i.e.}\
$h(x) = 1-x -\ldots - x^{k-1} + ax^k + x^{k+1}/(1-x)$ for some $k \ge 1$
and $a\in \{-1,0,1\}$. 
\end{theorem}

\noindent {\bf Remarks.} 1. Note that all the points described in the theorem are algebraic, since the function $h$ is rational over $\Z$.
The first ``tips of the corners'' to which the theorem applies, for an appropriate $(*)$-function, are as follows (given with 5-6 digit accuracy):

$(0.618034, 0.68232)$, which is a pair of zeros of $1 - x - x^2 - x^3 + x^5/(1 - x)$ (incidentally, the reciprocals of this pair are the golden ratio and the 4th Pisot number);

$(0.550607, 0.7691)$, which is a pair of zeros of $1 - x - x^2 - x^3 - x^4 + x^5/(1 - x)$;

$(0.532958, 0.804916)$, which is a pair of zeros of $1 - x - x^2 - x^3 - x^4 + x^6/(1-x)$;

$(0.519703, 0.83221)$, which is a pair of zeros of $1 - x - x^2 - x^3 - x^4 -x^5 + x^6/(1-x)$;

$(0.513951, 0.85068)$, which is a pair of zeros of $1 - x - x^2 - x^3 - x^4 -x^5 +x^7/(1-x)$.


2. The ``cusp corners'' obtained from $(*)$-functions are only the
``most outward'' cusp corners of $\Nk_+$. There are many others visible in Figures 1 and 2, which are probably
pairs of zeros of power series $h\in\Bk$
with all but finitely many  coefficients equal to $1$, as in Theorem~\ref{th-tips}.
For instance, it appears that there is a corner at $(0.645200, 0.68232)$ (with the second zero again the reciprocal of the 4th Pisot number), which is
a pair of zeros of $h(x)=1-x-x^2-x^3+x^4+x^6+x^8/(1-x)$. In order to prove this rigorously, one only needs to check that $h$ is the unique function in $\Bk$ with this pair of zeros, but we haven't done this.

\medskip

\begin{figure}[ht!]
\centering
\includegraphics[width=90mm]{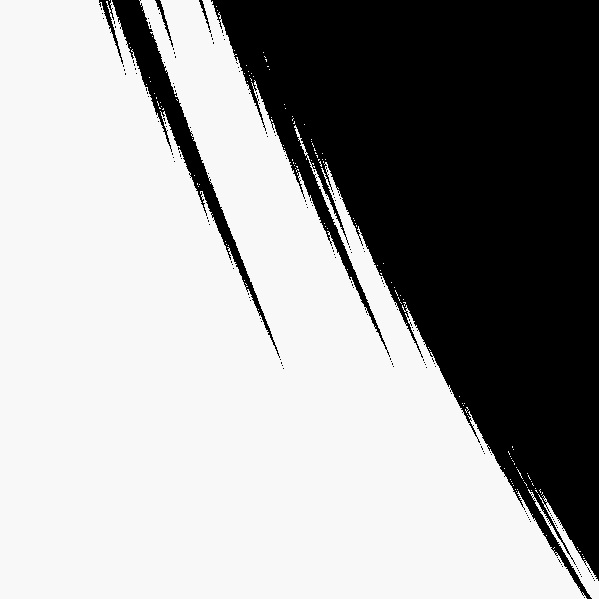}
\caption{The set $\Nk_+\cap [0.647,0.661]\times[0.677,0.691]$, with several prominent ``cusp corners''}
\label{overflow}
\end{figure}

\medskip

Our second main result is concerned with connectedness properties of the set $\Nk$.
Let 
$$
\Nktil_+ := \left\{(\gam,\lam) \in \Nk_+:\ \exists f \in \Bk,\
f(\gam) = f(\lam) =0,\ \xi_3(f) \le \lam \right\}.
$$
Recall that $\Nk_t=\{(\gam,\lam)\in (-1,1)^2: |\gam\lam|\ge \half\}$ is the ``trivial'' part of $\Nk$. .

\begin{theorem} \label{th-locon}
The set $\Nk_+ \setminus (\Nktil_+ \cup \Nk_t)$ is locally connected. Moreover,
there is no connected component of $\Nk_+$ that is disjoint from $\Nktil_+ \cup \Nk_t$.
\end{theorem}

We were  not able to prove the connectedness of the entire set $\Nk$, but conjecture that this is the case. 
The next proposition shows that the last theorem is non-vacuous,
in fact, $\Nk_+ \setminus (\Nktil_+ \cup \Nk_t)$ contains a substantial portion
of the set $\Nk_+\setminus \Nk_t$. In particular, the set $\Nk$ is connected near the ``cusp corners'' from Theorem~\ref{th-tips}.
Let
$$
\Omtil_+:= \left\{(\gam,\lam) \in \Om_+:\ \exists\,f \in \Bk_{[-1,1]},\
f(\gam) = f(\lam) = 0,\ \xi_3(f) \le \lam\right\}.
$$
Clearly, $\Nktil_+ \subset \Omtil_+$. Recall that $\alpha_3$ denotes
the smallest triple
zero of a power series in $\Bk_{[-1,1]}$; in \cite{BBBP} it is shown that
$\alpha_3 \approx .727883$ is a zero of a polynomial with integer
coefficients of degree 12.

\begin{prop}\label{prop-om2}
{\bf (a)} The function
$$
\psi:\ \gam\mapsto \min\{\xi_3(f):\ f \in \Bk_{[-1,1]},\ 
f(\gam) =0\}
$$
is well-defined on $(.5, \alpha_3)$. It is continuous,
decreasing, and satisfies
$$
\lim_{t\to .5+} \psi(t) = 1,\ \ \ \lim_{t\to \alpha_3-} \psi(t) = \alpha_3.
$$

{\bf (b)} For every $\gam \in (.5,\alpha_3)$ we have
$$
\gam < \lam < \psi(\gam)\ \ \Longrightarrow\ \ (\gam,\lam) \not\in \Omtil_+.
$$

{\bf (c)}  For every $\gam \in (.5,\alpha_3)$ there exists a unique
$(**)$-function
\be \label{eq-2star}
H_{k,\ell}^{(a,b)}(x) 
= 1 - \sum_{i=1}^{k-1} x^i + ax^k + \sum_{i=k+1}^{\ell-1} x^i  +
bx^\ell - \frac{x^{\ell+1}}{1-x} \in \Bk_{[-1,1]}
\ee
such  that $H_{k,\ell}^{(a,b)}(\gam) = H_{k,\ell}^{(a,b)}(\psi(\gam))=
(H_{k,\ell}^{(a,b)})'(\psi(\gam))=0$.
\end{prop}

The following table contains some values of the function $\psi$
(rounded-off in such a way that the actual values of $\psi$ are slightly
larger):

\bigskip

\begin{center}
\begin{tabular}{|c|c|c|c|c|c|c|c|c|c|c|c|c|c|c|c|}  \hline
$\gam$ & .53 & .55 & .57 & .59 & .61 & .63 & .65 & .67 & .69 & .71 & .7278 \\ \hline
$\psi(\gam)$ & .877 & .85 & .832 & .815 & .799 & .785 & .771 & .759 & .747 & .736 & .7278 \\ \hline
\end{tabular}
\end{center}


\bigskip


\begin{figure}[htb]
\centering
\includegraphics[width=0.8\textwidth]{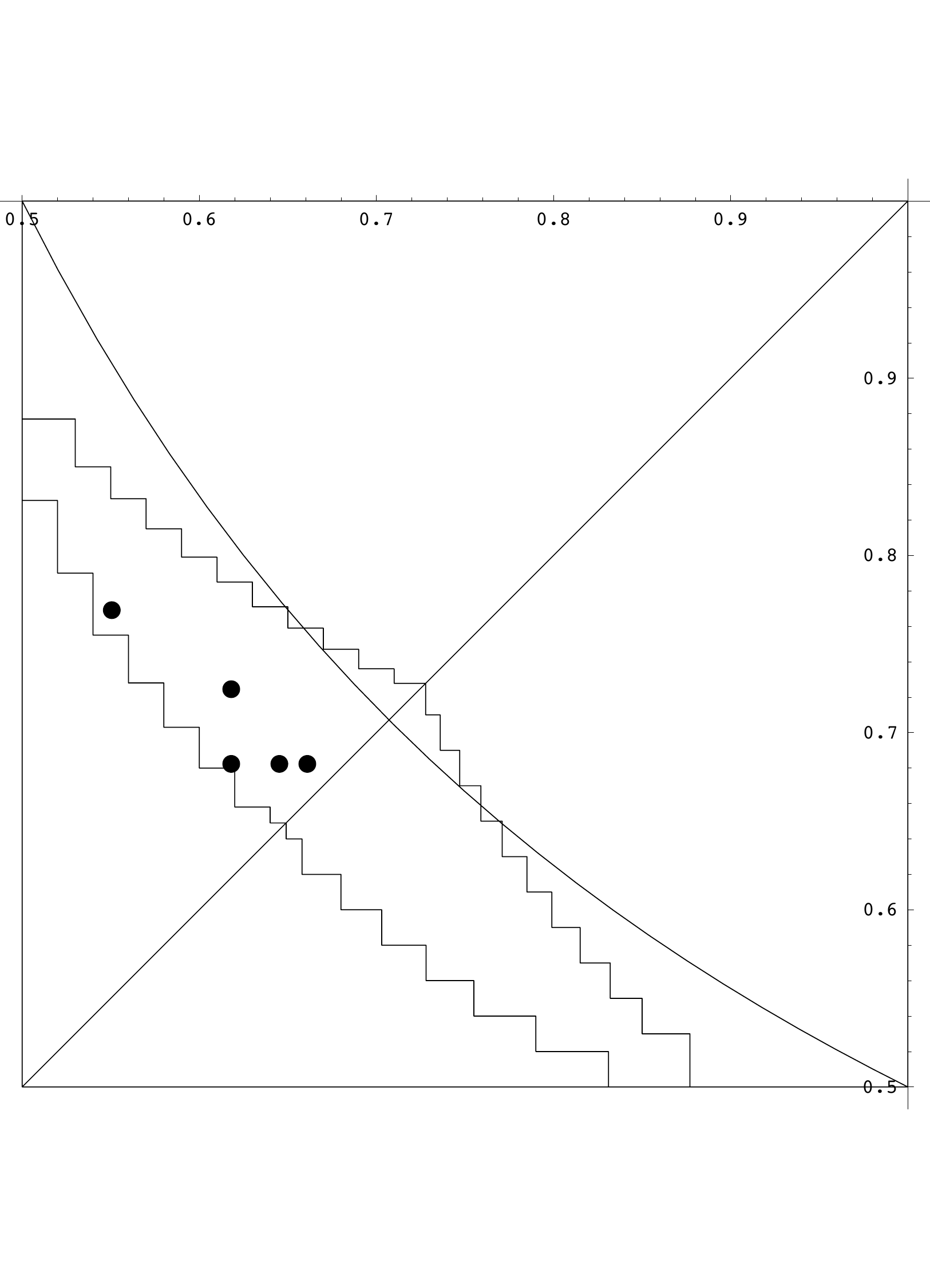}
\caption{}
\label{con-reg}
\end{figure}

Propositions~\ref{prop-om1} and \ref{prop-om2} are
illustrated in Figure 3, which shows the
region obtained from the tables. For $(\gam,\lam)$ in the region below
the lower broken line, the corresponding self-affine
set is totally disconnected. Theorem~\ref{th-locon} implies that the
part of $\Nk_+$ between the broken lines and outside $\Nk_t$ is locally connected and there are
no components of $\Nk$ entirely contained between the broken lines and
the diagonal. We also show several points which are known to belong to
$\Nk_+$ (these are some of the ``cusp corners'' from Theorem~\ref{th-tips}).

Another application of the bounds on the set 
$\Nktil_+$ comes from the paper by P. Shmerkin \cite{Pablo}. Following \cite[Def.\,4.10]{Pablo}, we say that $\Rk$ is a region of
$*$-transversality if for all $(\gam,\lam) \in \Rk$ there is $f\in \Bk$ such that $f(\gam) = f(\lam)=0$ but $f'(\gam) \ne 0,\ 
f'(\lam)\ne 0$. It is clear $\Nk_+\setminus \wtil{\Nk}_+$ is a region of $*$-transversality. It is proved in \cite{Pablo} that
for Lebesgue-a.e. $(\lam,\gam)$ with $\gam\lam < 1/2< \gam$ the self-affine attractor $K_{\gam,\lam}$ has Hausdorff dimension
$1 + \log(2\lam)/\log(1/\gam)$, and if $\Rk$ is a region of $*$-transversality contained in $\{(\gam,\lam) : \ \gam\lam < 1/2\}$, then 
$K_{\gam,\lam}$ has zero Hausdorff measure in its dimension for Lebesgue-a.e.\ $(\gam,\lam) \in \Rk$.


\section{Proofs}

\noindent {\em Proof of Proposition~\ref{prop-om1}.}
Consider the family of $(*)$-functions 
$h_k^{(a)}$ with $k \ge 1$ and $a \in [-1,1]$, given by (\ref{eq-star}),
and equip it with a total order as follows:
$$
h_k^{(u)} > h_\ell^{(v)}\ \ \mbox{if}\ \ k<\ell\ 
\ \mbox{or}\ \ k=\ell,\ u>v.
$$
Obviously,
$$
h_k^{(u)} > h_\ell^{(v)}\ \Longrightarrow\ h_k^{(u)}(x) > h_\ell^{(v)}(x)\ \
\mbox{for all}\ x\in (0,1).
$$
The set $\Bk_{[-1,1]}$ is a normal family of analytic functions in the
unit disk; therefore, it is compact in the uniform topology on any compact 
subset of $(0,1)$. We can also indentify $\Bk_{[-1,1]}$ with the
infinite product $[-1,1]^\infty$ equipped with the product topology.
Observe that $u\mapsto h_k^{(u)}$ is
a continuous function from $[-1,1]$ to $\Bk_{[-1,1]}$ for all $k\ge 1$.
It is strictly decreasing in the order defined above, and moreover,
$$
h_k^{(-1)} = h^{(1)}_{k+1}\ \ \ \mbox{for all}\ \ k\ge 1.
$$
By \cite{solerd} and \cite{BBBP}, there is a $(*)$-function
$h_4^{(b)}$ having a double zero at $\alpha_2$, with $b \approx .0875294$.

It is easy to see that every $h^{(u)}_k < h_4^{(b)}$ has exactly two
distinct positive zeros.
Indeed, $h^{(u)}_k$ has at least two positive zeros since $h^{(u)}_k(0)=1$,
$\lim_{x\to 1-} h^{(u)}_k(x) = +\infty$, and $h^{(u)}_k(\alpha_2)<
h_4^{(b)}(\alpha_2)=0$. On the other hand, the derivative
$(h^{(u)}_k)'(x)$ has only
one coefficient sign change, so it has at most one positive zero
by the Decartes Rule of Signs,
hence $h^{(u)}_k$ has at most two positive zeros.
Clearly, the zeros (when they exist) continuously depend on $u$.

\smallskip

{\sc Claim 1.} {\em For every $\gam \in (0.5,\alpha_2)$ there exists a
$(*)$-function $h^{(u)}_k$ such that $h^{(u)}_k(\gam)=0$.}

\smallskip

Indeed, $h_4^{(b)}(\gam)>0$, and for $k$ sufficiently large we have
$h^{(1)}_k(\gam)<0$, since $\lim_{k\to \infty} h^{(1)}_k(\gam) = 
1-\gam-\gam^2-\ldots = \frac{1-2\gam}{1-\gam}$.
By continuity, there exist $k$ and $u$ such that $h_k^{(u)}(\gam)=0$, and
of course, $h_k^{(u)} < h_4^{(b)}$.
Since $h_k^{(u)}(\alpha_2)<0$, there is another zero $\lam > \alpha_2>\gam$.
The claim is proved.

\smallskip

{\sc Claim 2.}  {\em We have $\phi(\gam) = \lam$, and $h_k^{(u)}$ is the unique
function in $\Bk_{[-1,1]}$ with zeros at $\gam$ and $\lam$.}

\smallskip

Indeed, suppose $f \in \Bk_{[-1,1]}$
is such that $f(\gam)=0$ and $f\not\equiv h_k^{(u)}$.
Consider $g(x)= f(x) - h_k^{(u)}(x)$. Then $g(x)$
is a power series with at most one coefficient sign change and $g(\gam)=0$.
It follows that $\gam$ is the only positive zero of $g$. The first nonzero
coefficient of $g$ is positive, so it is positive for small positive $x$.
It follows that $g(x)<0$ for all $x>\gamma$, hence 
$f(x) = g(x) + h_k^{(u)}(x)<0$ for all $x\in (\gam,\lam]$. So, for all $f \not\equiv h_k^{(u)} \in \Bk_{[-1,1]}$, $\xi_2(f) > \lam$, and the claim
is proved.

\smallskip

The claims show that the function $\phi$ is well-defined on $(\half,\alpha_2)$.
The remaining statements of part (a) are now easy to derive. In fact, one can obtain sharp asymptotics for $\phi(t)$ as $t\to \half+$ and
$t\to 1-$, but we do not pursue this.

\smallskip

(b) This statement is immediate from the definition of $\phi$.

\smallskip

(c) The formula for the ``optimal'' function and its uniqueness are already 
proved. 
The statement about the derivative of $h_k^{(a)}$ is also clear: as already mentioned, $(h_k^{(a)})'$
has only one sign change and its zero (which the minimum of $h_k^{(a)}$) must lie in $(\gam,\phi(\gam))$. \qed

\medskip

The table of $\phi$ values is obtained from Claims 1 and 2 above. See Appendix for details. 

\medskip

\noindent {\em Proof of Theorem~\ref{th-tips}.}
Suppose that $z_n=(\gam_n,\lam_n) \in \Nk_+$ are such that
$z_n \to z_0=(\gam_0,\lam_0)$. Consider functions (maybe non-unique)
$h_n \in \Bk$ such that
$h_n(\gam_n) = h_n(\lam_n)=0$. Since $\Bk$ is compact, there is 
a subsequence of $h_n$ converging to some $\htil \in \Bk$, with
$\htil(\gam_0)=\htil(\lam_0)=0$. By the assumption of uniqueness of such a function, we have $\htil = h$.
Since convergence in $\Bk$ is coefficientwise, it follows that for
any $N\in \Nat$ there is $n_0\in \Nat$ such that $h_n$ agrees with
$h$ in the first $N$ terms for all $n \ge n_0$. 

This already implies
that $z_0$ is a ``corner'' with interior angle at most $\pi/2$.
Indeed, if $h_n$ agrees with $h$
in the first $N\ge  N_0$ terms, where $N_0-1$ is the last term of $h$ with a coefficient different from $+1$, then $h-h_n$ has only non-negative coefficients
and hence $h_n(x)< h(x)$ for all $x\in (0,1)$. Since $h(x) \le 0$ for
$x\in [\gam_0,\lam_0]$ we obtain that the zeros of $h_n$ must satisfy
$\gam_n < \gam_0,\ \lam_n > \lam_0$.
For the more delicate estimate we need the following lemma:

\begin{lemma} \label{near-zeros} Suppose that $\half<\gam_0<\lam_0< 1$ are such that  $\gam_0$ and $\lam_0$ are zeros 
of  $h\in \Bk$, as in the statement of Theorem~\ref{th-tips}, and let $f(x) = h(x) - x^N R(x)$, where $R$ is a power series with 
coefficients $0,1$. Then for $N$ sufficiently large, $f$ has zeros $\gatil$ and $\latil$ satisfying
$$
\frac{2{\gam_0}^{N}R(\gam_0)}{3|h'(\gam_0)|} \leq \gam_0-\gatil  \leq \frac{2{\gam_0}^{N}R(\gam_0)}{|h'(\gam_0)|}\,,\ \ \ 
\frac{2{\lam_0}^{N}R(\lam_0)}{3h^{\prime}(\lam_0)} \leq \latil - \lam_0 \leq \frac{{2\lam_0}^{N}R(\lam_0)}{h^{\prime}(\lam_0)}\,.
$$ 
Recall that $h'(\gam_0)<0$ and $h'(\lam_0)>0$ by assumption.
\end{lemma}

First we deduce the theorem, assuming the lemma.
The argument at the beginning of the proof shows that for any $N\in \N$ there exists $\delta > 0$ such that for all $z = (\gam, \lam) \in \Nk_+ \cap B_{\delta}(z_0)$ (where $z_0 = (\gam_0, \lam_0)$), 
if $f \in \Bk$ is such that $f(\gam) = f(\lam) = 0$, then 
$f$ agrees with $h$ in the first $N$ coefficients. (Note that $z \in \Nk_+$ implies there does indeed exist $f \in \Bk$ such that 
$f(\gam)=f(\lam) = 0$.) 
Let $N_0\in \N$ be such that $h$ has only coefficients equal to $+1$ starting from $N_0$. Let $\delta>0$ be so small that $N\ge N_0$.
Then
$f(z) = h(z) - z^NR(z)$ for some power series $R$ with coefficients $0,1$ and we obtain from Lemma~\ref{near-zeros}:
$$
\frac{2{\gam_0}^{N}R(\gam_0)}{3|h^{\prime}(\gam_0)|} \leq \gam_0 - \gam \leq \frac{2{\gam_0}^{N}R(\gam_0)}{|h^{\prime}(\gam_0)|}\,, \ \ \ 
\frac{2{\lam_0}^{N}R(\lam_0)}{3h^{\prime}(\lam_0)} \leq \lam - \lam_0 \leq \frac{{2\lam_0}^{N}R(\lam_0)}{h^{\prime}(\lam_0)}\,.
$$ 
Let $\alpha = \frac{\log{\lam_0}}{\log{\gam_0}}$. Note that  ${\gam_0}^{\alpha} = \lam_0$, so
$$
\left(\frac{2}{3}\right)^{\alpha} \frac{{\lam_0}^NR(\gam_0)^{\alpha}}{|h^{\prime}(\gam_0)|^{\alpha}}  \leq  \left(\gam_0 - \gatil \right)^{\alpha} \leq \frac{2^{\alpha}{\lam_0}^NR(\gam_0)^{\alpha}}{|h^{\prime}(\gam_0)|^{\alpha}}\,.
$$
Thus we have that
$$
\frac{2R(\lam_0)|h^{\prime}(\gam_0)|^{\alpha}}{2^{\alpha}3R(\gam_0)^{\alpha}h^{\prime}(\lam_0)} \leq \frac{\latil - \lam_0}{\left(\gam_0 - \gatil \right)^{\alpha}} \leq \frac{3^{\alpha}2R(\lam_0)|h^{\prime}(\gam_0)|^{\alpha}}{2^{\alpha}R(\gam_0)^{\alpha}h^{\prime}(\lam_0)}\,.
$$
Now, $1 \leq R(\gam_0) \leq \frac{1}{1-\gam_0}$ and $1 \leq R(\lam_0) \leq \frac{1}{1-\lam_0}$, whence
$$
\frac{2|h^{\prime}(\gam_0)|^{\alpha}(1-\gam_0)^{\alpha}}{2^{\alpha}3h^{\prime}(\lam_0)} \leq \frac{\latil - \lam_0}{\left(\gam_0 - \gatil \right)^{\alpha}} \leq \frac{3^{\alpha}2|h^{\prime}(\gam_0)|^{\alpha}}{2^{\alpha}(1-\lam_0)h^{\prime}(\lam_0)}\,,
$$
as desired. The claim that the conditions on $h$ are satisfied whenever $h$ is a $(*)$-function is immediate from definitions and Proposition~\ref{prop-om1}. \qed

\medskip

{\em Proof of Lemma~\ref{near-zeros}.}
This is standard, but we provide the argument for completeness. We will only prove the estimate for $\latil$, since the one for $\gatil$ is obtained
in exactly the same way. 
We will need an easy inequality:
\be \label{easy}
|g''(x)| \le 2(1-x)^{-3}\ \ \ \mbox{for all}\ g\in \Bk\ \mbox{and}\ x\in (0,1).
\ee

 Recall that $h(\lam_0)=0$ and $h'(\lam_0)>0$, so $h(x)>0$ to the right of
$\lam_0$. Thus, it is clear that for large $N$ there will be a zero of $f(x) = h(x)-x^NR(x)$ in a small neighborhood 
$(\lam_0,\lam_0+t]$. Since the claim is local, we can assume that $\lam_0+t\le 1-\delta<1$ for some $\delta>0$ (independent of $N$, 
e.g.\ we can take $\delta = (1-\lam_0)/2$).

We have $f(\lam_0 + t) = h(\lam_0 + t) - (\lam_0 + t)^NR(\lam_0 + t)$. Recall that $f(\lam_0)<0$, and we want to make sure that
$f(\lam_0+t)\ge 0$. By Taylor's formula,
$$h(\lam_0 + t) \geq h^{\prime}(\lam_0)t - \frac{C_2t^2}{2},$$ where, in view of (\ref{easy}),
$$
C_2 :=2(1-\delta)^{-3}\ge \max\left\{|h''(x)|:\ x\in [\lam_0,\lam_0+t]\right\}.
$$
We can asume that $N$ is large enough, so that $$t := \frac{4\delta^{-1}\lam_0^N}{h'(\lam_0)} < \frac{4h'(\lam_0)}{C_2}\,.
$$
Then $h^{\prime}(\lam_0) t- \half C_2t^2 > \frac{1}{2}h^{\prime}(\lam_0)t$. We claim that  $ \frac{1}{2}h^{\prime}(\lam_0)t \geq (\lam_0 + t)^{N}R(\lam_0+t)$ for $N$ sufficiently large.
By the definition of $t$, 
\begin{eqnarray*}
(\lam_0 + t)^{N}R(\lam_0+t) & \le & \left(\lam_0 + \frac{4\delta^{-1}\lam_0^N}{h'(\lam_0)}\right)^N R(1-\delta) \\
& \le & \lam_0^N \left(1+ \frac{4\delta^{-1}\lam_0^{N-1}}{h'(\lam_0)}\right)^N \delta^{-1}.
\end{eqnarray*}
Since $\lim_{N\to \infty} \left(1+ 4\delta^{-1}\lam_0^{N-1}/h'(\lam_0)\right)^N =1$, we conclude that 
$(\lam_0 + t)^{N}R(\lam_0+t)\le 2\delta^{-1}\lam_0^N = \frac{1}{2}h^{\prime}(\lam_0)t$ for $N$ sufficiently large, as desired.
Thus we have shown that $f$ has a zero $\latil \in (\lam_0, \lam_0+t]$, where $t=4\delta^{-1}\lam_0^{N}/h'(\lam_0)$ and
$N$ is large enough.

	By the Mean Value Theorem, there exists $c \in (\lam_0, \latil)$ such that 
$
\lam_0^N R(\lam_0) = f(\latil) - f(\lam_0) = (\latil-\lam_0) f'(c),
$
so that
\be \label{eq-lati}
 \latil - \lam = \frac{\lam_0^N R(\lam_0)}{f'(c)}\,.
\ee
In view of (\ref{easy}),
$$|f'(c) - f'(\lam_0)| \le  |c-\lam_0| \cdot 2\delta^{-3} \\
 <  t \cdot 2\delta^{-3} = \frac{8\delta^{-4} \lam_0^N}{h'(\lam_0)}\,.
$$
Thus we can choose $N$ sufficiently large, so that
$|f^{\prime}(c) - f^{\prime}(\lam_0)| < \frac{1}{4}h^{\prime}(\lam_0)$. 
Next,
$$
|f'(\lam_0) - h'(\lam_0)| = |\lam_0^N R'(\lam_0) + N \lam_0^{N-1} R(\lam_0)|\le \lam_0^N \delta^{-2} + N \lam_0^{N-1} \delta^{-1},
$$
which is also less than $\frac{1}{4}h^{\prime}(\lam_0)$ for $N$ sufficiently large.
Then $|f^{\prime}(c) - h^{\prime}(\lam_0)| <  \frac{1}{2}h^{\prime}(\lam_0)$, whence
$f'(c) \in (\half h'(\lam_0), \frac {3}{2} h'(\lam_0))$, and so  we have
from (\ref{eq-lati}):
 $$
\frac{2{\lam_0}^{N}R(\lam_0)}{3h^{\prime}(\lam_0)} \leq \frac{{\lam_0}^{N}R(\lam_0)}{f^{\prime}(c)} = \latil - \lam_0 \leq \frac{{2\lam_0}^{N}R(\lam_0)}{h^{\prime}(\lam_0)}\,,
$$ 
as desired. \qed

\medskip

\noindent {\em Proof of Theorem~\ref{th-locon}.}
This proof is a modification of the argument by Bandt \cite[Section 11]{bandt}, which is, in turn, based on \cite{bou}.
Let $$\Delta:=\left\{(\gam,\lam) \in (0,1)^2:\ \gam<\lam,\ \gam\lam<1/2\right\}$$ and 
consider the quotient space
\be \label{defX}
\Xk:= \clos(\Delta\setminus \Nktil_+) / \partial(\Delta\setminus \Nktil_+)
\ee
with induced topology. Denote by $\om$ the point corresponding to the
contracted boundary.

Recall that $\Bk$ is the set of all power series of the
form $f(x) = 1+ \sum_{n=1}^\infty a_n x^n$, with $a_n \in \{-1,0,1\}$.
We can identify $\Bk$ with the space $\{-1,0,1\}^{\Nat}$ equipped
with the product topology. Observe that this topology coincides with the
topology of uniform convergence on compact subsets of the unit disk.

\smallskip

%
%
%

{\sc Claim 1.} {\em $\Nktil_+$ is relatively closed in $\Delta$.}

\smallskip

ndeed, let $(\gam_n,\lam_n) \to (\gam,\lam)\in \Delta$ and 
$(\gam_n,\lam_n)\in \Nktil_+$. Then
there exist $f_n \in \Bk$ with $f_n(\gam_n) = f_n(\lam_n) =0$ and
$\xi_3(f_n) \le \lam_n$. This means that there exist
$\alpha_n \le \lam_n$ such that $f_n(\alpha_n)=0$
(if $\alpha_n$ is equal to $\gam_n$ or $\lam_n$ this is understood as
having the corresponding zero of multiplicity 2).
By compactness, without loss of generality, we
can assume that $f_n \to f\in \Bk$ and $\alpha_n \to \alpha$. Then
$f(\gam) = f(\lam) = f(\alpha)=0$ and $\alpha \le \lam$
(again using our convention concerning double zeros). Thus $(\gam, \lam) \in \Nktil_+$, and Claim 1 is proved.

\smallskip

We will also need the following fact (see \cite{BBBP} and \cite[Th.\,2.4]{Pablo}): if $f\in \Bk$ and $\alpha_1,\ldots,\alpha_k$ are (some) complex roots of $f$ in the unit disk, counted with multiplicity, then
$
|\alpha_1\cdots \alpha_k| \ge (1+k^{-1})^{-k/2} (k+1)^{-1/2}.
$
Taking $k=4$, we obtain
\be \label{rootb}
|\alpha_1\alpha_2\alpha_3 \alpha_4|\ge 16\cdot 5^{-5/2}> 1/4.
\ee

Let $\phi:\ \Bk \to \Xk$ be the function defined as follows:
If $f \in \Bk$ is such that $\gam = \xi_1(f) < \xi_2(f)=\lam$ and
$(\gam,\lam) \in \Delta \setminus \Nktil_+$, then
$\phi(f) := (\gam,\lam)$; otherwise, $\phi(f) := \om$.

\smallskip

{\sc Claim 2.} {\em $\phi:\ \Bk \to \Xk$ is continuous.}

\smallskip

Indeed, if $\phi(f) = (\gam,\lam)$,
then $f$ has simple zeros at $\gam$ and $\lam$, hence a small perturbation
of $f$ will result in a small perturbation of these zeros.
Suppose that $\phi(f) = \om$. 
We need to show that if $\ftil$ is
a small perturbation of $f$, then $\phi(\ftil)$ is close to $\om$.
We have the following possibilities:

\smallskip

(a) $f$ has no positive zeros;

(b) $f$ has one simple positive zero;

(c) $\xi_1(f) < \xi_2(f) < \xi_3(f)$ but
$(\xi_1(f),\xi_2(f)) \in \Nktil_+$;

(d) $\xi_1(f)=\xi_2(f)< \xi_3(f)$; 

(e) $\xi_1(f) < \xi_2(f) = \xi_3(f)$;

(f) $\xi_1(f) = \xi_2(f) = \xi_3(f)$, that is, $f$ has a triple zero.

\smallskip

It is not hard to see that for any $\eps>0$ there exists $\delta>0$ such
that if the distance from $\ftil$ to $f$ in $\Bk$ is less than $\delta$,
then every zero of $\ftil$ in $(0,1-2\eps)$ is $\eps$-close to a 
zero of $f$. Another general fact is useful: for a small perturbation of a real-analytic function, new real zeros cannot appear; zeros can only
dissappear (i.e.\ become non-real). 

In cases (a) and (b), either $\ftil$ has the same
property, hence $\phi(\ftil)=\om$, or new zeros appear near 1, which
could result in $\phi(\ftil) = (\gam,\lam)$ with $\lam$ near 1, that is,
$(\gam,\lam)$ is close to $\om$ in the topology of $\Xk$.
In the case (c), a small perturbation $\ftil$ will have
$\phi(\ftil)=\om$ or $\phi(\ftil)=
(\xi_1(\ftil),\xi_2(\ftil))$, which is
close to $(\xi_1(f),\xi_2(f)) \in \Nktil_+$,
hence close to $\om$ in the topology of $\Xk$. Suppose that case (d) holds. If the double zero at $\xi_1(f) = \xi_2(f)$ doesn't dissappear, then
we either still have a double zero for $\ftil$, or two real zeros close to each other. In the former case we have $\phi(\ftil) = \om$, and in the
latter case $\phi(\ftil)=
(\xi_1(\ftil),\xi_2(\ftil))$ is close to the diagonal, that is, close to $\om$ in the topology of $\Xk$. If, on the other hand, the double zero dissappears (becomes non-real) and 
$\phi(\ftil)\ne\om$, then $\phi(\ftil) = (\xi_1(\ftil),\xi_2(\ftil))$, which is close to 
$(\gam,\lam)$ where $\gam \le \lam$ are zeros of $f$ and
$\xi_1(f)= \xi_2(f) < \gam$. However, in the latter case we have $\gam\lam>1/2$ by (\ref{rootb}), hence
$(\gam,\lam) \not \in \Delta$ and so $(\xi_1(\ftil), \xi_2(\ftil))$ is still close to $\om$. If (e) or (f) holds, we have
a similar argument: assuming $\phi(\ftil)=(\gatil,\latil)\ne\om$, this point is either close to one
in $\Nktil$, or close to the diagonal, or close to $\partial \Delta$, in view of (\ref{rootb}). This concludes the proof of Claim 2.

\smallskip

Now we essentially repeat the argument from \cite{bou}. For a finite
word $u = u_1\ldots u_n$ in the alphabet $\{-1,0,1\}$ let
$$
\Bk_u = \left\{ 1 + u_1 x + \ldots + u_n x^n + \sum_{k=n+1}^\infty b_k x^k \in
\Bk\right\}.
$$
Denote by $u,\alpha$ the word of length $n+1$ obtained by adding the
symbol $\alpha$ to $u$.
Let $f = 1 + u_1 x + \ldots + u_n x^n \in \Bk_{u,0}$. We have
$f(1-x^{n+1})^{-1} \in \Bk_{u,1},\ f(1-x^{n+1})^{-1} \in \Bk_{u,-1}$.
Observe that these three functions have the same set of zeros in $(0,1)$,
hence they are mapped into the same point by $\phi$. It follows that
$$
\phi(\Bk_{u,-1}) \cap \phi(\Bk_{u,0}) \cap \phi(\Bk_{u,1}) \ne \es
$$
for an arbitrary $u$.
This property is called {\em recursive connectedness} in \cite{bou}.
It is proved in \cite{bou} (and quite easy to see) that this property
implies that $\phi(\Bk)$ is connected and locally connected in $\Xk$.
Note that $\phi(\Bk) = \pi(\Delta\cap\Nk_+ \setminus \Nktil_+) \cup \{\om\}$
where $\pi:\ \clos(\Delta\cap \Nk_+ \setminus \Nktil_+) \to \Xk$ is the natural
projection associated with the quotient map. (We know that
$\om\in\phi(\Bk)$ since, for instance,
there are power series in $\Bk$ with no positive zeros.)
This immediately implies the statement of the theorem. \qed

\medskip

{\em Proof of Proposition~\ref{prop-om2}.}
By \cite[Cor.\,2.5.]{Pablo-Sol}, for any $\gam>\half$ there exists $f\in \Bk
\subset \Bk_{[-1,1]}$ such that $f(\gam)=0$ and $f$ has at least three
zeros in $(0,1)$ (for example, we may take the zeros of $f$ to be $\gam$ with multiplicity one, and $\sqrt{\frac{1}{2\gam}}$ with multiplicity two). Thus the set in the definition of the function $\psi$ 
is non-empty; it has a minimum by compactness of the class $\Bk_{[-1,1]}$.
The statement (b) is immediate from the definitions.
The remaining statements of (a) will easily follow from part (c). 
Its proof is divided into several lemmas.

\begin{lemma} \label{lem-doubs}
Suppose that $h$ is a $(**)$-function such that $h'(x_0) =0$
and $h(t)>0,\ h'(t)<0$ for all $t \in (0,x_0)$. Then there is no $f\in
\Bk_{[-1,1]}$ such that $\xi_3(f)=x_0$, unless $h=f$, and
$x_0$ is a triple zero of $f=h$.
\end{lemma}

\noindent {\em Proof.} Suppose that $f\in \Bk_{[-1,1]}$ violates the
assertion of the lemma. 
Let $g(x) = f(x)-h(x)$. By the definition of a $(**)$-function, we have
$$
g(x) = A_1(x) - A_2(x) + A_3(x),
$$
where $A_1(x)$ and $A_2(x)$ are polynomials and $A_3(x)$ is a power series,
all three
with non-negative coefficients, such that the highest power in 
$A_i$ is less than the lowest power in $A_{i+1}$. Thus, $g$ and $g'$
have at most two coefficient sign changes each.

Since $x_0$ is the third zero of $f$ and $f(0)=1$, we have $f'(x_0)\le 0$. (Indeed, otherwise $f$ is negative in a left neighborhood of $x_0$, but 
on an interval where a real-analytic function changes its sign it must have an odd number of zeros, counting with multiplicities.)
Thus, $g'(x_0)=f'(x_0) - h'(x_0) = f'(x_0)\le 0$. Observe that there must be a zero $\zeta_1$
of $f'$ between the
first and second zeros of $f$ (if these two zeros of $f$ coincide, that is,
it is a double zero, which is equal to $\zeta_1$). We have $f(\zeta_1)\le 0$,
hence $g(\zeta_1) = f(\zeta_1) - h(\zeta_1)<0$, and
$g'(\zeta_1) = -h'(\zeta_1)>0$. By the Decartes Rule of Signs, $g'$ can have at most
two positive zeros. There has to be a zero of $g'$ in $(\zeta_1,x_0)$. 
There also have to be another zero of $g'$ in $(0,\zeta_1)$, since $g(0)=0$, $g(\zeta_1)<0$, and $g'(\zeta_1)>0$.
 Thus, $g'$ has exactly two
coefficient sign changes, hence $A_1(x)\not\equiv 0$. 
But then $g$ increases sufficiently close to zero, whence $g'$ must have
at least {\em two} zeros in $(0,\zeta_1)$. This is a contradiction. \qed

\begin{lemma} \label{lem-alpha3} $\alpha_3 = \min\{\xi_3(f):\
f\in \Bk_{[-1,1]}\}$.
\end{lemma}

\noindent {\em Proof.} Suppose there exists $f \in \Bk_{[-1,1]}$ such that
$\lam:=\xi_3(f) < \alpha_3$. 
It is proved in \cite{BBBP} that there is a $(**)$-function
$H=H_{k,\ell}^{(u,v)}$ such that $H(\alpha_3) = H'(\alpha_3)
= H''(\alpha_3)=0$. (In fact, $k=4$ and $\ell=10$.)   Consider the function
$$
h(x)  = H(x) + sx^\ell,\ \ \ \mbox{with}\ \ s = 
-\frac{H'(\lam)}{\ell\lam^{\ell-1}}\,.
$$
This is a $(**)$-function, though not necessarily in $\Bk_{[-1,1]}$, since the $x^\ell$-coefficient may exceed 1 in absolute value.
We have $H(x)>0,\ H'(x)<0$ for $x\in (0,\alpha_3)$, hence
$h(x) > 0$ for all $x\in (0,\lam)$ and $h'(\lam)=0$ by definition.
We claim that $h'(x)<0$ for all $x\in (0,\lam)$. Indeed, $h'$ has
two coefficient sign changes, hence at most two positive zeros.
We know that $h'$ is negative near zero, $h'(\lam)=0$, $h'(\alpha_3)=s\ell\alpha_3^{\ell-1}>0$, and
$h'$ is negative sufficiently close to 1. It follows that $h'$ has  a zero
in $(\alpha_3,1)$, so it does not vanish in $(0,\lam)$, implying the claim.
Thus, $h$ is a
$(**)$-function satisfying the assumptions of Lemma~\ref{lem-doubs}
for $x_0=\lam$, so the existence of $f$ is a contradiction. \qed

\medskip

Lemma~\ref{lem-alpha3}  implies that $\psi(\gam)\ge \alpha_3>\gam$ for 
$\gam \in (\half,\alpha_3)$.
Fix $\gam \in (\half,\alpha_3)$.
Recall that $\psi(\gam)$ is well-defined, which means that there exists
a function $f\in \Bk_{[-1,1]}$ such that $\psi(\gam) = \xi_3(f)$. Such
a function will be called ``optimal'' (for a given $\gam$).

\smallskip

\begin{lemma} \label{lem-claim1} An optimal function $f$ for $\gam\in (\half,\alpha_3)$ has a double zero at
$\lam = \xi_3(f)$, that is, $f(\lam) = f'(\lam)=0$.
\end{lemma}

\noindent {\em Proof.}
Suppose $f'(\lam)\ne 0$. Since $f(0)=1$ and $\lam$ is the third
positive zero of a real analytic function, $f$ is strictly decreasing
in a neighborhood of $\lam$. By Decartes' Rule of Signs, $f$ has at least
three coefficient sign changes. Therefore, we can find  integers
$0 < \ell_1 < \ell_2 < \ell_3$ such that $a_{\ell_1} < 0,\
a_{\ell_2} > 0,$ and $a_{\ell_3} < 0$, where $a_{\ell_i}$ is the coefficient of $x^{\ell_i}$ in $f$.
Consider
$$
\ftil(x):= f(x) + \eps (\gam^{\ell_2-\ell_1} x^{\ell_1} - x^{\ell_2}).
$$
Then $\ftil \in \Bk_{[-1,1]}$ for sufficiently small $\eps>0$. 
Moreover, $\ftil(\gam)=f(\gam)=0$ and  $\ftil(x)<f(x)$ for $x\in (\gam,1)$.
Thus, for sufficiently small $\eps>0$, the function
$\ftil$ has a zero close to $\lam$ which is less than $\lam$. We claim that
this zero is $\xi_3(\ftil)$, which contradicts $\lam = \psi(\gam)$.
Indeed, if the first two positive zeros of $f$ are distinct (and they are
smaller than $\lam$), this property will persist for $\ftil$. If $\gam$ is
a double zero, then $\ftil$ has a second zero $\gam'$ close to $\gam$.
This proves the claim, and the lemma follows. \qed

\begin{lemma} \label{lem-claim2}
The optimal function $f$ for $\gam\in (\half,\alpha_3)$ is unique; it is a
$(**)$-function $h_{k,\ell}^{(a,b)}$ for some 
$1 \le k < \ell < \infty$ and $a,b \in [-1,1]$.
\end{lemma}

\noindent {\em Proof.}
Let $f(x) = 1+\sum_{n=1}^\infty a_n x^n$ 
be optimal, and suppose that it is not a $(**)$-function 
(see (\ref{eq-2star})).
Let $\ell_1\ge 1$
be minimal such that $a_{\ell_1} > -1$. Then choose $\ell_2>\ell_1$
minimal such that $a_{\ell_2} < 1$ (note that $\ell_2$ exists since $f$
must have at least three coefficient sign changes). 
If $f$ is not a $(**)$-function, then we can find $\ell_3>\ell_2$ such that
$a_{\ell_3} > -1$. Let $c_2,c_3\in \R$ be such that 
$$
g(x) := -x^{\ell_1} + c_2 x^{\ell_2} + c_3 x^{\ell_3}
$$
satisfies $g(\gam) = g(\lam)=0$. (This is a linear system of equations with determinant $\gam^{\ell_2}\lam^{\ell_3} - \gam^{\ell_3}\lam^{\ell_2}\ne 0$,
so there is a unique solution.)
Notice that $c_2>0$ and $c_3<0$, since there must be two coefficient sign changes in
$g$. Clearly, $\lam$ is a simple zero for $g$,  and it is a double zero for
$f$ by Lemma~\ref{lem-claim1}. Thus, there exist $b_1,b_2>0$ such that
$$
|f(x)| \le b_1|x-\lam|^2,\ \ \ \ \ |g(x)| \ge b_2|x-\lam|
$$
for $x$ near $\lam$. Consider
$$
\ftil(x):= f(x) + \eps g(x).
$$ 
Then $\ftil \in \Bk_{[-1,1]}$ for sufficiently small $\eps>0$. Observe that
$\ftil(\gam) = \ftil(\lam)=0$ by construction. Recall that $f(0)=1,
f(\gam)=0$, and $f(\lam)= f'(\lam)=0$, hence $\min_{[\gam,\lam]} f < 0$. We
can make sure that $\eps>0$ is so small that 
$\min_{[\gam,\lam]} \ftil <0$. On the other hand,
$$
\ftil(\lam-\eps/n) = f(\lam-\eps/n) + \eps g(\lam-\eps/n) \ge -b_1(\eps/n)^2 + \eps b_2 (\eps/n)>0
$$
provided that $n > b_1/b_2$. Then $\ftil$ has a zero in $(\gam,\lam)$ which
implies that $\xi_3(\ftil) \le \lam$. Since $\lam = \psi(\gam)$, we have
$\xi_3(\ftil)= \lam$, so $\ftil$ is optimal for $\gam$ (as well as $f$). 
This contradicts Lemma~\ref{lem-claim2} since $\lam$ is not a double zero
of $\ftil$.  It remains to verify that the optimal function is unique.
Assuming that we have two distinct optimal functions, we take their difference,
which has at most two 
coefficient sign changes, since both are $(**)$-functions.
This leads to a contradiction since the difference has at least three
positive zeros. 

This concludes the proof of the lemma and of the claim (c) in Proposition~\ref{prop-om2}. The remaining statements of the proposition
follow easily. \qed 


\section{Appendix: how to compute the functions $\phi$ and $\psi$}

We first explain how the function $\phi$ was computed, using Mathematica.
Consider the $(*)$-function $h(x) = 1 - x - . . . - x^{k-1} + ax^k + \frac{x^{k+1}}{1-x}$, where $a \in [-1,1]$. 
First, we fix $k$. The algorithm takes $\gam$ as an input. Then, $a = F(\gam)$ is determined so that $h(\gam) = 0$. We must check that $-1 \leq a \leq 1$, so that $h$ is indeed a member of $\Bk_{[-1,1]}$. Next, we find the second root of $h$ using the FindRoot command with an appropriate starting point. We choose the starting point to guarantee that we find $\lam$ rather than $\gam$ (FindRoot uses Newton's method to find the root of a function. It will find the root closest to the starting point. Recall the shape of the $(*)$-function $h$. We must choose a starting point to the right of the minimum of $h$ to guarantee that Mathematica finds $\lam$ rather than $\gam$. We know our choice of starting point works as long as the output of FindRoot is not equal to our input $\gam$).

\begin{figure}[htb] \label{h_1}
\centering
\includegraphics[width=0.6\textwidth]{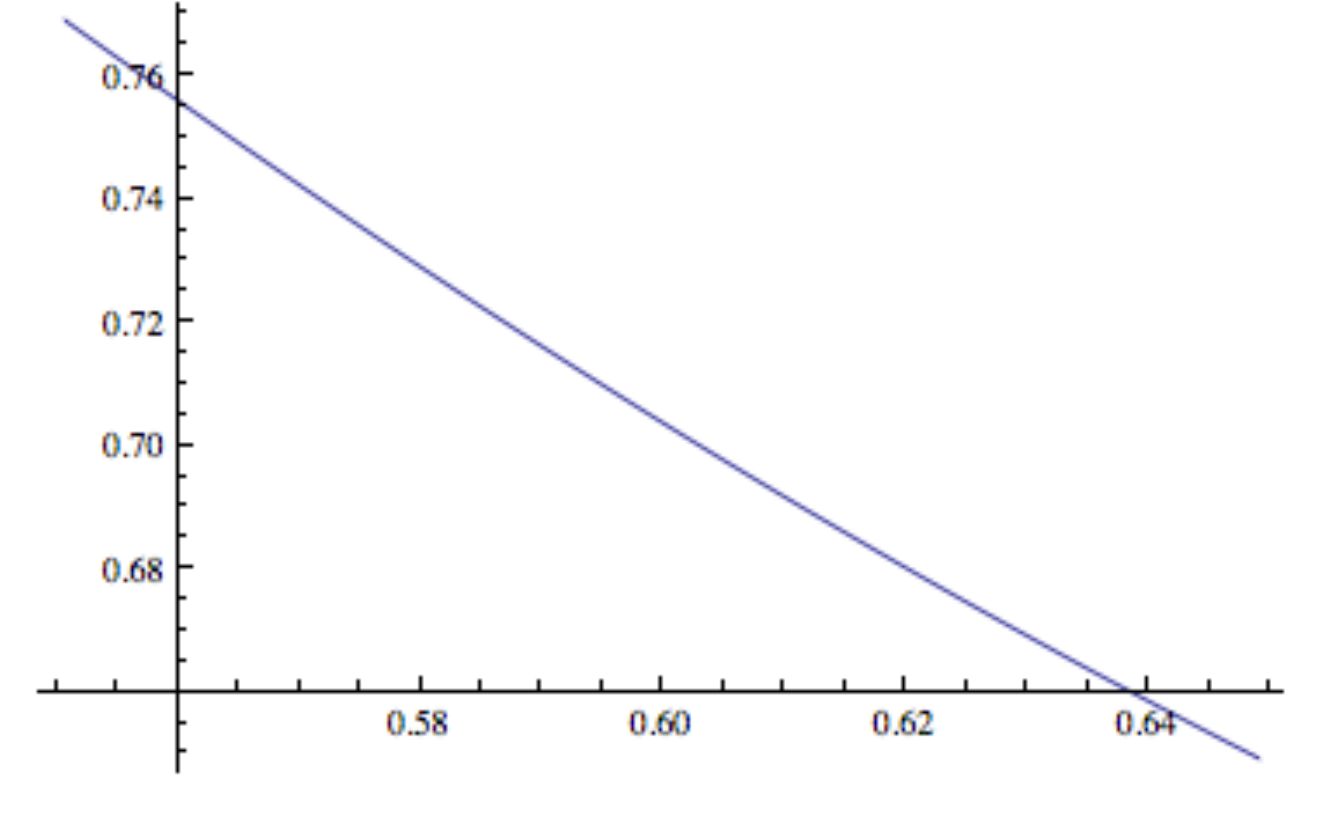}
\caption{$\phi(\gam)$ from $h_1$}
\end{figure}

 Now, as was seen in \cite{solerd,BBBP}, 
there is a $(*)$-function $h_4^{(b)}$ having a double zero at $\alpha_2$. Therefore we begin by fixing $k = 4$. 
Consider the $(*)$-function $h_1(x)  = 1 - x - x^2 - x^3 - F(\gam)x^4 + \frac{x^5}{1-x}$. We solve for $F(\gam)$ to ensure that $h_1(\gam)=0$
and obtain
$F(\gam) = \frac{1-2\gam+\gam^4 + \gam^5}{\gam^4-\gam^5}$. 
We find that 0.8 works as a starting point for FindRoot.
Using NSolve, we find that $|a| = |F(\gam)| \leq 1$ for $\gam \in (0.550607, 0.7691)$, approximately. However, recall that $\phi : \, (0.5, \alpha_2) \rightarrow (0,1)$. Thus we are only interested in looking at $\gam \in (0.550607, \alpha_2 = 0.649138)$.  Figure \ref{h_1} shows a plot of $\phi(\gam)$ for $\gam \in (0.550607, 0.649138)$.

Now, note that when $\gam = 0.550607$, $F(\gam) \approx 1$.  Thus at $\gam = 0.550607$, the coefficient of $x^4$ is $-1$. Thus this is a ``switching point," that is, at this point, $h_1$ switches to a $(*)$-function with $k = 5$. This is one of the points we are interested in, because it will be in the set $\Nk$. 

Next we consider the $(*)$-function $h_2(x) = 1 - x - x^2 - x^3 - x^4 - G(\gam)x^5 + \frac{x^6}{1-x}$, so that $k = 5$. We solve for $G(\gam)$ so that $\gam$ is indeed a root of $h_2(x)$, and find $G(\gam) = \frac{1-2\gam + \gam^5 + \gam^6}{\gam^5 - \gam^6}$.
We again check the range for which $|G(\gam)| \leq 1$, and find that the inequality holds for $\gam \in (0.519703, 0.832218)$.
\begin{figure}[htb] \label{h_2}
\centering
\includegraphics[width=0.6\textwidth]{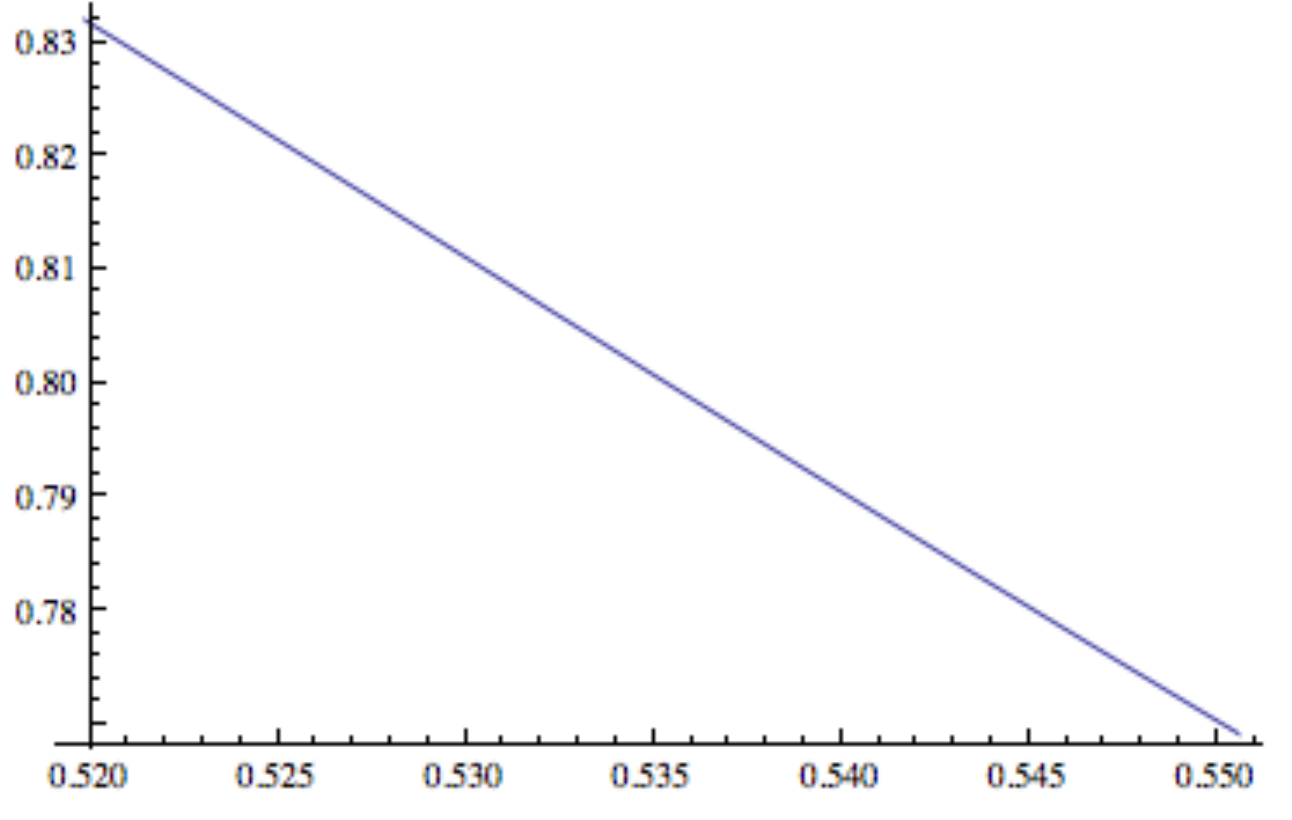}
\caption{$\phi(\gam)$ from $h_2$}
\end{figure}
\begin{sloppypar}
Figure 5 shows a plot of $\phi(\gam)$ for $\gam \in (0.519703, 0.550607)$.
Note that $G(0.529703) \approx 1$, so that at $\gam = 0.529703$ $h_2$ becomes a $(*)$-function with $k = 6$. Thus we continue similarly by setting $k = 6$. Let $h_3(x) = 1-x-x^2-x^3-x^4-x^5 - K(\gam)x^6 + \frac{x^7}{1-x}$, and solve for $K(\gam)$ so that $h_3(\gam) = 0$ to obtain
$K(\gam) = \frac{1-2\gam +\gam^6 + \gam^7}{\gam^6 - \gam^7}$.
Using NSolve, we find that $|K(\gam)| \leq 1$ for $\gam \in (0.508831, 0.866368)$. 
We may continue in this manner in order to obtain $\phi(\gam)$ for $\gam \rightarrow 0.5$. Note, however, that the process does not terminate.
\end{sloppypar}

\medskip

Now we explain how the function $\psi$ may be computed.
Recall that $\psi : (0.5, \alpha_3) \rightarrow [0,1]$, where $\alpha_3 \approx 0.727883$. 
Consider the $(**)$-function $$\displaystyle H_{k,\ell}(x) = 1 - \sum_{i=1}^{k-1}x^i + ax^k + \sum_{i = k+1}^{\ell-1}x^i + bx^\ell- \frac{x^{\ell+1}}{1-x}\,,$$
where $a,b \in [-1,1]$. 
Recall that we would like to find $H_{k,\ell}(x)$ such that $H_{k,\ell}(\gam) = H_{k,\ell}(\lam) = H_{k,\ell}^{\prime}(\lam) = 0$. Since we have two unknowns $a$ and $b$, for this algorithm we will start with $\lambda$ and obtain $a = F_a(\lam)$ and $b = F_b(\lam)$ such that $H_{k,\ell}(\lam) = H_{k,\ell}^{\prime}(\lam) = 0$, and use FindRoot to find $\gam$ such that $H_{k,\ell}(\gam) = 0$. 
In \cite{BBBP} it was proved that there is a $(**)$-function H=$H_{4,10}$ such that $H(\alpha_3) = H^{\prime}(\alpha_3) = H^{\prime \prime}(\alpha_3) = 0$.
 Thus we begin by considering the function $F(x) = 1-x-x^2-x^3+x^5+x^6+x^7+x^8+x^9 - \frac{x^{11}}{1-x}$. Then we  let $H_1(x) = F(x) + F_a(\lam)x^4 + F_b(\lam)x^{10}$.
 Note that $F^{\prime}(x) = -1 - 2x - 3x^2 + 5x^4 + 6x^5 + 7x^6 + 8x^7 + 9x^8 - \frac{11x^{10}}{1-x} - \frac{x^{11}}{(1-x)^2}$. 
We solve the system of equations
\begin{eqnarray*}
 H_1(\lam) &=& F(\lam) + F_a(\lam)\lam^4 + F_b(\lam)\lam^{10} = 0\\
 H_1^{\prime}(\lam) &=& F^{\prime}(\lam) + 4F_a(\lam)\lam^3 + 10F_b(\lam)\lam^9 = 0
\end{eqnarray*}
 for $F_a(\lam)$ and $F_b(\lam)$ and find that $F_a(\lam) = \frac{\lam F^{\prime}(\lam)-10F(\lam)}{6\lam^4}$ and $F_b(\lam) = \frac{4F(\lam)-\lam F^{\prime}(\lam)}{6\lam^{10}}$. 
In this case, we need to have both $|F_a(\lam)| \leq 1$ and $|F_b(\lam)| \leq 1$.
We find that $|F_a(\lam)| \leq 1$ for $\lam \in (0.606471, 0.83611)$, (note that $F_a(0.606471) = F_a(0.83611) = -1$) and $|F_b(\lam)| \leq 1$ for $\lam \in (0.692945, \alpha_3)$ (where $F_b(0.692945) = -1$). So, the first coeffcient that ``switches" is at $k = 10$, when $F_b(\lam) = -1$. 
Next we use FindRoot to find $\gamma$. So this algorithm takes $\lambda$ as an input and outputs $\gamma$, so this is effectively $\psi^{-1}$. 

Next, we use our ``switching point." We let $F_1(x) = F(x) + x^{10} + x^{11}$ and $G_1(x) = F_1(x) + F_{a_1}(\lam)x^4 + F_{b_1}(\lam)x^{11}$ and  solve for $F_{a_1}(\lam)$ and $F_{b_1}(\lam)$ so that $G_1(\lam) = G^{\prime}_1(\lam) = 0$. We find that $F_{a_1}(\lam) = \frac{\lam F_1^{\prime}(\lam) - 11F_1(\lam)}{7\lam^4}$ and $F_{b_1}(\lam) = \frac{-\lam F_1^{\prime}(\lam) + 4F_1(\lam)}{7\lam^{11}}$. 
Again we check for which $\lam$ we have $|F_{a_1}(\lam)| \leq 1$ and $|F_{b_1}(\lam)| \leq 1$ and continue in this way.

%

\medskip

\noindent
{\bf Acknowledgment.}
This paper was started in 2003, but took a long time complete. Its preliminary drafts are cited in \cite{Pablo,Pablo-Sol}. Thanks to Pablo Shmerkin for helpful discussions and fruitful collaboration on \cite{Pablo-Sol}. I am grateful to Christoph Bandt for his help with the computer program
(which he specially wrote on my request), without which this paper would have been impossible. Thanks also to Allison
Beckwith who helped  with the Appendix, as part of her Master's project at the University of Washington Department
of Mathematics during 2011-2012. 


\bibliographystyle{amsplain}

\end{document}